\documentclass[12pt]{article}

\usepackage[latin2]{inputenc}
\usepackage{epsfig}
\usepackage{latexsym}
\usepackage{graphicx}
\usepackage{makeidx}
\usepackage{amsmath,amssymb,amsbsy,amsfonts,amsthm,latexsym,amsopn,amstext,amsxtra,euscript,amscd}
\usepackage{xcolor}
\usepackage{arydshln}

\newtheorem{theorem}{Theorem}

\newtheorem{corollary}{Corollary}

\newtheorem{remark}{Remark}

\def\bg{\bar{g}}
\def\br{\bar{r}}
\def\G{\mathbf{G}}
\def\M{\mathbf{M}}
\def\T{\mathbf{T}}
\def\M{\mathbf{M}}

\def\B{\mathbf{B}}
\def\A{\mathbf{A}}

\def\O{\mathbf{O}}

\def\N{\mathbb{N}}
\def\Z{\mathbb{Z}}
\def\Q{\mathbb{Q}}
\def\R{\mathbb{R}}
\def\C{\mathbb{C}}

\begin{document}

\centerline{{\bf Decomposable forms generated by linear recurrences}}

\centerline{K. Gy\H ory\footnote{University of Debrecen (Hungary)}, A. Peth\H o\footnote{University of Debrecen (Hungary)}, L. Szalay\footnote{University of Sopron (Hungary), and J.~Selye University (Slovakia)}}
\vspace{1cm}

\begin{abstract}
Consider $k\ge 2$ distinct, linearly independent, homogeneous linear recurrences of order $k$ satisfying the same recurrence relation. We prove that the recurrences are related to a decomposable form of degree $k$, and there is a very broad general identity with a suitable exponential expression depending on the recurrences. This identity is a common and wide generalization of several known identities. Further, if the recurrences are integer sequences, then the diophantine equation associated to the decomposable form and the exponential term has infinitely many integer solutions generated by the terms of the recurrences. We describe a method for the complete factorization of the decomposable form. Both the form and its decomposition are explicitly given if $k=2$, and we present a typical example for $k=3$. The basic tool we use is the matrix method.
\end{abstract}



\section{Introduction}\label{Introduction}

Let $k\ge2$ be an integer. A sequence $(G_n)$ of complex numbers  is called {\it linear recurrence of order $k$} if there exist $\gamma_0,\dots,\gamma_{k-1}\in \mathbb{C}$ such that
\begin{equation} \label{e:rek}
G_{n+k} = \gamma_{k-1} G_{n+k-1} + \dots + \gamma_0G_n
\end{equation}
holds for all $n\ge 0$, moreover $k$ is minimal with this property.  The polynomial 
\begin{equation}\label{cp}
	P(X) = X^k -\gamma_{k-1} X^{k-1}-\dots-\gamma_0
\end{equation} 
is the {\it characteristic polynomial} of $(G_n)$. Note that the definition of the order yields $\gamma_0\not= 0$. We will return to this question later, in Section 2.

Dividing (\ref{e:rek}) by $\gamma_0\not=0$ and rearranging the equality we obtain
$$
G_n = -\frac{\gamma_1}{\gamma_0}G_{n+1}-\dots -\frac{\gamma_{k-1}}{\gamma_0}G_{n+k-1}+\frac{1}{\gamma_0}G_{n+k}.
$$
Now performing the substitution $n \mapsto -n-k$ it gives
$$
G_{-(n+k)} = -\frac{\gamma_1}{\gamma_0}G_{-(n+k-1)}-\dots -\frac{\gamma_{k-1}}{\gamma_0}G_{-(n+1)}+\frac{1}{\gamma_0}G_{-n},
$$
which is a linear recursion for the negative branch of $(G_n)$. The results of the present paper hold for this two side infinite sequence what we can consider also as a number theoretic function $G_n \,:\, \mathbb{Z} \mapsto \mathbb{C}$.

There exist a huge number of relations involving linear recurrences in the corresponding literature, in particular for $k=2$. Identities with Fibonacci $(F_n)$, or Lucas $(L_n)$ numbers play distinguished role among them. The Fibonacci Quarterly is devoted mainly for such works. Sequences $(F_n), (L_n)$ have the initial terms $F_0=0, L_0=2, F_1=L_1=1$, and both satisfy the recursion
\begin{equation}\label{x}
x_n=Ax_{n-1}+Bx_{n-2}
\end{equation}
with $A=B=1$. From this rich variety of identities we quote three ones, namely
\begin{eqnarray} 
L_n^2 - 5 F_n^2 &=& (L_n - \sqrt{5}F_n)(L_n + \sqrt{5}F_n) = 4(-1)^n,  \label{e:Fib1}\\
F_{n+1}^2-F_nF_{n+1}-F_{n}^2 &=& \left(F_{n+1} - \frac{1+\sqrt{5}}{2}F_n \right) \left(F_{n+1} - \frac{1-\sqrt{5}}{2}F_n\right) = (-1)^n, \label{e:Fib2}\\
F_{n-1}F_{n+1}-F_n^2&=& (-1)^n. \label{e:Fib3}
\end{eqnarray}
All three belong to the folklore and can be easily generalized to any second order recurrences. Indeed, assume that $(G_n)$, and $(\widehat{G}_n)$ satisfy the recursion \eqref{x}. If their initial terms are $G_0,G_1, \widehat{G}_0=2G_1-AG_0,\widehat{G}_1=AG_1+2BG_0$, then we say that $(\widehat{G}_n)$ is {\it associated} to $(G_n)$. Here $A,B, G_0,G_1$ denote complex numbers with the conditions $B\not= 0$ and $|G_0|+|G_1|\not= 0$. Assume that 
$\alpha$ and $\beta$ are the not necessarily different roots of the characteristic polynomial $P(X)=X^2-AX-B$. Put $D=A^2+4B$, which is the discriminant of $P$, and let $C_G=G_1^2-AG_0G_1-BG_0^2$. Then the identities
\begin{eqnarray} 
\widehat{G}_n^2-DG_n^2&=&(\widehat{G}_n-\sqrt{D}G_n)(\widehat{G}_n+\sqrt{D}G_n) = 4C_G(-B)^n, \label{Gen1}\\
G_{n+1}^2-AG_nG_{n+1}-BG_{n}^2 &=& \left(G_{n+1} - \alpha G_n \right) \left(G_{n+1} - \beta G_n\right) = C_G(-B)^{n}, \label{e:Gen2}\\
G_{n-1}G_{n+1}-G_n^2&=&-C_G(-B)^{n-1} \label{jan13}
\end{eqnarray}
hold, and it is a simple exercise to see that they are direct generalizations of \eqref{e:Fib1}-\eqref{e:Fib3}, respectively. Note that the second one is equivalent to the third one by $(G_{n+1}-AG_n)G_{n+1}-BG_n^2=BG_{n-1}G_{n+1}-BG_n^2$, and this phenomenon is obviously true for \eqref{e:Fib2} and \eqref{e:Fib3}, too.

Assuming $|B|=1$ and using \eqref{Gen1}, Peth\H{o} \cite{Petho} characterized the polynomials whose values appear infinitely many times in $(G_n)$. One can reverse  identity \eqref{e:Fib2} such that the only positive integer solutions to the equation 
$$
x^2-xy-y^2 = \pm1
$$
are $(F_n,F_{n+1})$. J. P. Jones \cite{Jones} gave a good overview on this equation and its relation to the solution of Hilbert's tenth problem.

For recurrences of higher order there are only a few and complicated identities  available. Especially, we do not know any generalization of \eqref{Gen1}. Recently, Craveiro et al.~\cite{CSR} explored a nice generalization of the so-called Cassini identity (\ref{e:Fib3}) for recurrences of arbitrary order $k$. They investigated also analytical and combinatorial properties of their result. Corollary \ref{c:ugyanaz} of Theorem \ref{t:fo} of this paper provides a new contribution to this question, in fact we consider the problem with any set of initial values apart from singular cases.  Corollary \ref{c:ugyanaz} via Theorem \ref{t:fo} makes it possible to determine a Cassini-like identity if the basic recurrence of order $k$ is given. However such an identity is expected to be rather compound.

The main result of this article is a common and wide extension of \eqref{Gen1}, \eqref{e:Gen2}, and other identities to  $k\ge2$ linear recurrences of order $k$ satisfying the same recursive rule, including their connection with decomposable forms and diophantine equations.
 
\section{New results}

We shall prove some general results. Note that the basic field we work in is the set $\mathbb{C}$ of complex numbers, but the machinery works even for any fields.

\begin{theorem} \label{t:fo}
Let $(G_n^{(j)}), j= 1, \dots,k$ be linear recursive sequences of complex numbers of order $k$ with the same characteristic polynomial (\ref{cp}) having constant term $-\gamma_0\not= 0$. If the sequences are $\mathbb{C}$-linearly independent, then there exists a homogeneous polynomial $F\in \mathbb{C}[X_1,X_2,\dots,X_k]$ of degree $k$ such that 
$$
F(G_n^{(1)}, G_n^{(2)},\dots, G_n^{(k)}) = ((-1)^{k+1}\gamma_0)^n
$$
holds for all $n\in \mathbb{Z}$. 
\end{theorem}

A homogeneous polynomial $Q\in \mathbb{C}[X_1,X_2,\dots,X_k]$ is called {\it decomposable} if it splits completely into linear factors. Hence a decomposable form is necessarily homogeneous, but for $k\ge 3$ the converse is not true in general. The next theorem is a remarkable completion of Theorem \ref{t:fo}.

\begin{theorem} \label{t:decompose}
Under the assumptions of Theorem \ref{t:fo} the homogeneous polynomial $F$ is decomposable.
\end{theorem}

Clearly, Theorem \ref{t:decompose} is a continuation of Theorem \ref{t:fo}, and we could have presented them together. However their proofs differ significantly, which explains their separation.  

We emphasize that the proofs are constructive, and following their arguments one can compute the polynomial ${F}$ and its decomposition into linear factors, too. Unfortunately, $F$ has no nice general form if $k\ge3$, so we will omit its explicit constitution. On the other hand, $F$ is given precisely when $k=2$ (see Theorem \ref{main2} and the description after). Furthermore, we will show some representative examples in Section \ref{s:negy}.

A simple consequence of the main theorems is the following generalization of \eqref{e:Gen2}.

\begin{corollary} \label{c:ugyanaz}
Let $(G_n)$ be a linear recursive sequence of complex numbers of order $k$ satisfying \eqref{e:rek}. If the vectors $\bar{g}_j=(G_j,G_{j+1},\dots,G_{j+k-1}), j=0,\dots,k-1$ are $\mathbb{C}$-linearly independent, then there exists a decomposable form $F\in \mathbb{C}[X_1,X_2,\dots,X_k]$ of degree $k$ such that 
$$
F(G_n,G_{n+1}, \dots, G_{n+k-1}) = ((-1)^{k+1}\gamma_0)^n
$$
holds for all $n\in \mathbb{Z}$. 
\end{corollary}

The linear recursive sequences with initial values $G_0=\dots=G_{k-2}=0, G_{k-1}=1$ satisfy always the assumption of this corollary, hence the statement too. The $k$-generalized Fibonacci sequences  $(F_n^{(k)})$ (or, in other words, $k$-step Fibonacci numbers), which are defined by the initial values $F_0^{(k)} = \dots = F_{k-2}^{(k)} = 0, F_{k-1}^{(k)} =1$ and by the recursion $F_{n+k}^{(k)}=F_{n+k-1}^{(k)}+\dots+F_n^{(k)}$ are important examples. Here the upper index $(k)$ in $F_n^{(k)}$  means traditionally the parameter $k$, while in Theorem \ref{t:fo} the upper index $(j)$ denotes the $j$th from the $k$ given sequences. Clearly, the case $k=2$ is the Fibonacci sequence.

In fact, sequence $(G_n)$ is a complex valued number theoretical function. Another direct consequence of 
Theorem \ref{t:fo} is as follows. The conditions of Theorem \ref{t:fo} are valid here, too.

\begin{corollary} \label{c:indep}
Let $(G_n^{(j)}), j= 1, \dots,k$ be linear recursive sequences of complex numbers of order $k$ with the same characteristic polynomial (\ref{cp}). Then they and the sequence $(\gamma_0^n)$ are algebraically dependent.
\end{corollary}

It is possible to extend $(G_n)$ to a complex function $G(z)$ such that $G_n = G(n)$ for all $n\in \mathbb{Z}$. Indeed, assume that $P(X)$ in (\ref{cp}) has the factorization
$$
P(X) = (X-\alpha_1)^{m_1}(X-\alpha_2)^{m_2}\cdots (X-\alpha_h)^{m_h},
$$
where $\alpha_1,\dots,\alpha_h$ denote pairwise distinct complex numbers and $m_j$, $j=1,\dots,h$ are positive integers. Then, see e.g. \cite{ShT}, there exist polynomials $g_j\in \mathbb{C}[X]$ of degree at most $m_j-1$ such that
$$
G_n = g_1(n)\alpha_1^n +g_2(n)\alpha_2^n + \dots + g_h(n)\alpha_h^n
$$
holds for all $n\in \mathbb{Z}$. 
Hence $G(z)=g_1(z)\alpha_1^z +g_2(z)\alpha_2^z + \dots + g_h(z)\alpha_h^z$ is a complex function with $G(n) = G_n$ for all $n\in \mathbb{Z}$. For example, the extensions of the Fibonacci and the Lucas sequences are the complex functions  
$$
F(z) = \frac{1}{\sqrt{5}}\left(\left(\frac{1+\sqrt{5}}{2}\right)^z - \left(\frac{1-\sqrt{5}}{2}\right)^z\right),\quad L(z) = \left(\frac{1+\sqrt{5}}{2}\right)^z + \left(\frac{1-\sqrt{5}}{2}\right)^z,
$$
respectively. One can easily verify that identities \eqref{e:Fib1}-\eqref{e:Fib3} remain true if we replace $n$, $F_n$, $L_n$, $(-1)^n$ by $z$, $F(z)$, $L(z)$, $(-1)^z$, respectively. An interpretation of the generalized identities is the fact that the corresponding complex functions are algebraically dependent.

It is well known, see e.g. \cite{Konig}, that a function $f(z)$ satisfies the linear differential equation with characteristic polynomial $P(X)$ if and only if it is identical with one of the above defined $G(z)$. Our next theorem is a generalization of Corollary \ref{c:indep}.

\begin{theorem} \label{t:indep}
Let $P(X)\in \mathbb{C}[X]$ be a polynomial of degree $k$ with simple roots, and assume $P(0) \not= 0$. Let $k_0\ge k$ and denote $G_j(z)$, $1\le j\le k_0$ pairwise different solutions of a homogeneous linear differential equation with characteristic polynomial $P(X)$. Then  $G_j(z)$, $1\le j\le k_0$ and the function $((-1)^k P(0))^z$ are algebraically dependent.
\end{theorem}


As we noted Theorem \ref{t:decompose} and its constructive proof have a remarkable consequence in the theory of diophantine equations. This is presented as Theorem \ref{main} in Section \ref{Homogeneousform}, where it is more suitable to formulate the precise statement. Indeed, if $(G_n^{(j)})$'s are integer sequences, then it makes possible to have infinitely many integer solutions to the polynomial-exponential diophantine equation
\begin{equation} \label{e:normform}
\tilde{F}(x_1,x_2,\dots,x_k)=pq^n
\end{equation}
in the integers $x_1,x_2,\dots,x_k$ and $n\ge0$, where $\tilde{F}$ is a decomposable form of degree $k$ depending on the parameters of the sequences, further $p$, $q$ are suitable integers (see (\ref{theeq}) below). If $\tilde{F}$ is irreducible over $\mathbb{Q}[X_1,\dots,X_k]$, then we may assume that $\tilde{F}$ is, up to a constant factor, a full norm
form. It follows from a result of Borevich and Safarevich \cite{BorSaf} that if the splitting field of $\tilde{F}$ differs from $\mathbb{Q}$ and from the imaginary quadratic number fields, then \eqref{e:normform} has infinitely many solutions for $q=1$ and for some integer $p$. More generally, by a celebrated theorem of W. M. Schmidt \cite{Schmidt} on norm form equations the set of solutions is a union of a finite set and finitely many families of solutions. For a further generalization for arbitrary decomposable form equations see Gy\H{o}ry \cite{Gyory}. 
The families mentioned are sums of product of powers of appropriate algebraic numbers. Fixing all but one exponent and let run the remaining exponents through non-negative integers we see that each family includes linear recursive sequences. 

Our result shows that \eqref{e:normform} may have infinitely many solutions in the case when $\tilde{F}$ is reducible over $\mathbb{Q}[X_1,\dots,X_k]$, too. 

In our paper, we apply linear algebraic approach. We were also motivated to contribute to the development of so-called matrix method often used in the theory of linear recurrences. At the early 80's H. W. Gould \cite{G} presented a survey  on the $Q$-matrices, and he refers to R. Simson who first gave the formula
\eqref{e:Fib3}. This identity, together with its generalization \eqref{jan13} is an easy consequence of our Theorem \ref{main2} in Section 6 with the notation $H_n=G_{n+1}$ (for $H_n$ see again Theorem \ref{main2}). Theorem \ref{main2} is a common extension of identities (\ref{Gen1})-(\ref{jan13}), and provides the background for the case of $k\ge2$ linear recurrences of order $k$ in general appearing in Theorem \ref{t:fo}. 

We mention that recently Craviero et al. \cite{CSR} have given a generalization of Cassini identity for the $k$-generalized Fibonacci numbers with specific initial vector set. The basic approach of \cite{CSR} coincides with the idea of Lemma 2.1 of \cite{NSz}, and we exploit the advantages of this argument in the present paper.

The paper is organized as follows. In Section \ref{Homogeneousform}, we give two proofs for Theorem \ref{t:fo}. In the next section, we prove Theorem \ref{t:decompose}, while Section \ref{s:corrolaries} is devoted to the proof of corollaries \ref{c:ugyanaz} and \ref{c:indep} and of Theorem \ref{t:indep}. Section \ref{s:negy} specializes the results of Theorem \ref{t:fo} and \ref{t:decompose} to the cases $k=2$ and $k=3$. In Subsection 6.1, we present explicit formula in full generality if $k=2$, while in Subsection 6.2, we compute the appropriate formula only for three given third order recursive sequence. Other examples for larger $k$ values might be easily presented by following the method of the proofs.

\section{The homogeneous polynomials}\label{Homogeneousform}

Before starting the first proof of Theorem \ref{t:fo} we introduce some notation and summarize well known facts.

Recall that $k\ge2$ is an integer, furthermore $\gamma_0,\gamma_1,\dots,\gamma_{k-1}$ denote arbitrary complex numbers. Consider the set of recurrent sequences
\begin{equation}\label{Gamk}
\Gamma_k^{(\gamma_0,\gamma_1,\dots,\gamma_{k-1})}=\left\{(x_n)_{n\in\Z}\in\C^\infty\mid x_n=\gamma_{k-1}x_{n-1}+\gamma_{k-2}x_{n-2}+\cdots+\gamma_0x_{n-k}\right\}.
\end{equation}
Clearly, sequences (\ref{Gamk}) satisfy a common recurrence rule, but they differ in their initial values. They constitute a $\C$-vector space with respect to the coordinate-wise addition and multiplication by scalar if and only if $\gamma_0\not=0$, and in this case the dimension of the vector space is $k$. To be able to apply linear algebraic tools we assume $\gamma_0\not=0$ in the sequel. Throughout this paper we deal  with vector spaces over the field of complex numbers $\C$, therefore we do not mention always the ground field.

\subsection{First proof of Theorem \ref{t:fo}}

Assume that the recurrences
\begin{equation}\label{kseq}
(G_n^{(1)}),(G_n^{(2)}),\dots,(G_n^{(k)})\in\Gamma_k^{(\gamma_0,\gamma_1,\dots,\gamma_{k-1})}
\end{equation}
are $\mathbb{C}$-linearly independent. Then, equivalently, the vectors
\begin{equation}\label{defg}
\bg_0=\left[\begin{array}{c}
G_0^{(1)} \\ \vdots \\ G_0^{(k)} \\
\end{array}
\right],\;
\bg_1=\left[\begin{array}{c}
G_1^{(1)} \\ \vdots \\ G_1^{(k)} \\
\end{array}
\right],\dots,\;
\bg_{k-1}=\left[\begin{array}{c}
G_{k-1}^{(1)} \\ \vdots \\ G_{k-1}^{(k)} \\
\end{array}
\right]
\end{equation}
are also linearly independent. Define the matrices $\G,\G^\star\in\C^{k\times k}$ as follows. The column vectors of $\G$ are $\bg_0,\bg_1,\dots,\bg_{k-1}$, respectively, while $\bg_1,\bg_2,\dots,\bg_{k}$ admit the column vectors of $\G^\star$ in this order. The $k$ recurrences in (\ref{kseq}), together with their initial values build up the vector recurrence 
\begin{equation}\label{vecrec}
	\bg_n=\gamma_{k-1}\bg_{n-1}+\gamma_{k-2}\bg_{n-2}+\cdots+\gamma_0\bg_{n-k},\quad n\ge k
\end{equation}
with initial vectors (\ref{defg}).
In particular, (\ref{vecrec}) provides the last column vector $\bg_{k}$ of $\G^\star$.
The linear independence of the column vectors of $\G$ guarantees that
\begin{equation}\label{Delta}
	\Delta=\det(\G)\ne0.
\end{equation}

Let $\br_i$ denote the $i$th column vector of the transposition $\G^\top$ of $\G$ for each $i=1,2,\dots,k$. That is $\br_i=[G_0^{(i)},G_1^{(i)},\dots,G_{k-1}^{(i)}]^\top$, the entries are exactly the initial values of the $i$th sequence we fixed in (\ref{kseq}).
Clearly, the vectors
\begin{equation}\label{rbasis}
\br_1,\br_1\dots,\br_{k}
\end{equation}
form a basis of the vector space $\C^k$.

Put
\begin{equation*}
\M=\G^\star \G^{-1}=\frac{1}{\Delta}\cdot\G^\star\cdot{\rm adj}(\G).
\end{equation*}
The extension  of Lemma 2.1 of \cite{NSz} from $\R$ to $\C$ shows that the vector sequence $(\bg_n)$ generated by the initial vectors \eqref{defg} and by the vector recurrence (\ref{vecrec})
satisfies
\begin{equation} \label{NSz}
\bg_{n+1}=\M\bg_n,\qquad n\ge0.
\end{equation}
The proof of this lemma provides even that the characteristic polynomial of matrix $\M$ is
\begin{equation}\label{charpoly}
k_\M(x)=\det(x\mathbf{I}-\M)=x^k-\gamma_{k-1}x^{k-1}-\cdots-\gamma_1x-\gamma_0,
\end{equation}
which is also the characteristic polynomial of any recurrence of $\Gamma_k^{(\gamma_0,\gamma_1,\dots,\gamma_{k-1})}$ ($\mathbf{I}$ denotes the $k\times k$ unit matrix). Thus $k_\M(0)=\det(-\M)=-\gamma_0$, and then $\det(\M)=(-1)^k(-\gamma_0)$, is denoted in the sequel by $\delta$. Hence
\begin{equation}\label{detMn}
\det(\M^n)=\left(\det(\M)\right)^n=\delta^n.
\end{equation}
This observation gives later the right-hand side of the principal identity of this article (with a suitable homogeneous polynomial of degree $k$ on the left-hand side).
\smallskip

Now turn our attention to the matrices
\begin{equation*}
\M^n=\left[m_{i,j}^{(n)}\right]_{i,j=1,2,\dots,k},\qquad n\in\N.
\end{equation*}
Since $\M^n=\gamma_{k-1}\M^{n-1}+\gamma_{k-2}\M^{n-2}+\cdots+\gamma_0\M^{n-k}$ holds for $n\ge k$, therefore the same recurrence rule is valid for each element-wise sequence
$(m_{i,j}^{(n)})$ of the matrix sequence $(\M^n)$. That is
\begin{equation*}
m_{i,j}^{(n)}=\gamma_{k-1}m_{i,j}^{(n-1)}+\gamma_{k-2}m_{i,j}^{(n-2)}+\cdots+\gamma_0m_{i,j}^{(n-k)}.
\end{equation*}
Such a sequence is determined by the initial values $m_{i,j}^{(0)},m_{i,j}^{(1)},\dots,m_{i,j}^{(k-1)}$. Collect them into the initial column vector $\overline{m}_{i,j}\in\C^k$ (look at \eqref{mij}), which can be given as a linear combination of the basis vectors \eqref{rbasis}. More precisely, there exist uniquely determined coordinates $c_{i,j}^{(u)}\in\C$, $u=1,2,\dots,k$, such that
\begin{equation}\label{mij}
\overline{m}_{i,j}=\left[\begin{array}{c}
m_{i,j}^{(0)} \\
m_{i,j}^{(1)} \\ \vdots \\ m_{i,j}^{(k-1)}\\
\end{array}
\right]=\sum_{u=1}^{k} c_{i,j}^{(u)}\br_u.
\end{equation}
This is a system of $k$ linear equations in the $k$ unknowns $c_{i,j}^{(1)},c_{i,j}^{(2)},\dots,c_{i,j}^{(k)}$ (if $i,j$ are fixed), and the system can be solved, for instance, by Cramer's rule. Assume that $\G_u$ is the matrix derived by replacing the $u$th column vector of $\G^\top$ by the vector $\overline{m}_{i,j}$. Then,
using $\det(\M^\top)=\det(\M)=\Delta$,
clearly
\begin{equation}\label{coeffff}
c_{i,j}^{(u)}=\frac{\det(\G_u)}{\Delta}.
\end{equation}
Since the initial values determine a complete sequence of $\Gamma_k^{(\gamma_0,\gamma_1,\dots,\gamma_{k-1})}$, the coefficients (\ref{coeffff}) appearing in (\ref{mij}) descend for the whole sequence $(m_{i,j}^{(n)})$. Thus
the main consequence of the previous arguments is the equality
\begin{equation*}
m_{i,j}^{(n)}=\sum_{u=1}^{k} c_{i,j}^{(u)}G_n^{(u)}.
\end{equation*}

In this manner, we are able to represent any entry of matrix $\M^n$ as a linear combination of the $n$th terms of the linear recurrences
$(G_n^{(1)}),(G_n^{(2)}),\dots,(G_n^{(k)})\in\Gamma_k^{(\gamma_0,\gamma_1,\dots,\gamma_{k-1})}$. Subsequently, the determinant of $\M^n$ is a homogeneous polynomial of degree $k$ in $G_n^{(1)},G_n^{(2)},\dots,G_n^{(k)}$. Denote this form by $F$. Combining it with \eqref{detMn} we obtain
\begin{equation}\label{mainresult}
F(G_n^{(1)},G_n^{(2)},\dots,G_n^{(k)})=\delta^n,
\end{equation}
and the proof is complete.
$\Box$

Assume now that the coefficients $\gamma_i$, and the initial terms of the sequences $G_n^{(j)}$ are integers. Then $(G_n^{(j)})$ is a sequence of integers ($j$ is fixed). Multiply both sides of \eqref{mainresult} by $\Delta^k$, because (\ref{coeffff}) is now a rational number with denominator dividing $\Delta$.
Then $\tilde{F}=\Delta^k F$ has integer coefficients. Consequently, the polynomial-exponential diophantine equation
\begin{equation}\label{gede}
\tilde{F}(x_1,x_2,\dots,x_k)=\Delta^k\delta^n
\end{equation}
possesses infinitely many integer solutions in integers $x_1,x_2,\dots, x_k$ and $n$. Now we record the results in 

\begin{theorem}\label{main}
Given the recursive sequences (\ref{kseq}) of integers with initial values (\ref{defg}) such that (\ref{Delta}) holds. 
Then for any $n\ge0$
$$
(x_1,x_2,\dots,x_k)=(G_n^{(1)},G_n^{(2)},\dots,G_n^{(k)})
$$
is the solution to the diophantine equation
\begin{equation}\label{theeq}
\tilde{F}(x_1,x_2,\dots,x_k)=\Delta^k\delta^n,
\end{equation}
where the coefficients of the homogeneous form $\tilde{F}$ depend on the parameters of the sequences.
\end{theorem}

\begin{proof}
	The proof is given before the enunciation of the theorem.
\end{proof}

Note that the proof is constructive in the sense that following the method we can compute $\tilde{F}$, $\Delta$, and $\delta$ if the initial values and the recurrence rule is fixed. This will be illustrated later, in Section \ref{s:negy} for binary recurrences in general, and for a triple of ternary recurrences.

\subsection{Alternative proof of Theorem \ref{t:fo}}

We provide here an alternative proof for the existence of the form  ${F}$ on the left-hand side of (\ref{mainresult}). The notation is taken from the preceding parts of Section \ref{Homogeneousform}. Observe first that 
\begin{equation*} \label{rekurzio}
	{\bf G}^* = {\bf G} {\bf T}
\end{equation*}
with the matrix
$$
{\T} = \left[
\begin{array}{ccccc}
	0 & 0  & \cdots & 0 & \gamma_0 \\
	1 & 0  & \cdots & 0 & \gamma_1 \\
	\vdots & \vdots & \ddots & \vdots & \vdots \\
	0 & 0 & \cdots & 0 & \gamma_{k-2} \\
	0 & 0 & \cdots & 1 & \gamma_{k-1}
\end{array}
\right].
$$
Thus we have
$$
{\bf M} = {\bf G}^* {\bf G}^{-1} = {\bf G} {\bf T} {\bf G}^{-1}.
$$
This matrix equality immediately implies
\begin{itemize}
	\item $\det({\bf M}) = \det({\bf T}) = (-1)^{k-1}\gamma_0$, and
	\item ${\bf M}^n = {\bf G} {\bf T}^n {\bf G}^{-1}$ for $n\ge 0$.
\end{itemize}

Here  ${\bf G} {\bf T}^n$ is a matrix with the column vectors $\bar g_n, \bar g_{n+1},\ldots,\bar g_{n+k-1}$. On the other hand, from \eqref{NSz} we know $\bar{g}_{n+j} ={\bf M}^j \bar{g}_n$. Hence all entries of ${\bf G} {\bf T}^n$ are linear combinations of the terms $G_n^{(1)},G_n^{(2)},\ldots,G_n^{(k)}$, where the coefficients are the elements of ${\bf M}^j, j=0,\ldots,k-1$, and they are independent from $n$. Subsequently, the determinant of ${\bf G} {\bf T}^n$ is a signed sum of the products of $k$ linear forms in $G_n^{(1)},G_n^{(2)},\ldots,G_n^{(k)}$ i.e. a homogeneous form of degree $k$. Multiplying this by $\det({\bf G}^{-1})$, which is the determinant of a constant matrix we obtain that $\det({\bf M}^n)$ is also a homogeneous form of degree $k$. $\Box$

\section{Decomposability of $\mathbf{F}$}\label{Decomposability}

This section is devoted to study, and prove the decomposability of the homogeneous polynomial $F$ provided in Section \ref{Homogeneousform}. Recall that $F(x_1,x_2,\dots,x_k)=\tilde{F}(x_1,x_2,\dots,x_k)/\Delta^k$ (see (\ref{mainresult}) and (\ref{gede})).

\medskip
{\it Proof of Theorem \ref{t:decompose}.} 
The proof consists of two main parts. In the first step, we show that a specific homogeneous polynomial, say $F_S$ is decomposable if the initial $k$ sequences from $\Gamma_k^{(\gamma_0,\gamma_1,\dots,\gamma_{k-1})}$ are chosen advantageously. Then, in the second part, we find a connection between these favorable sequences and arbitrary package of $k$ sequences.
\bigskip

\noindent{\bf Part 1.}
Recall that the characteristic polynomial $k_\M(x)$ of the family $\Gamma_k^{(\gamma_0,\gamma_1,\dots,\gamma_{k-1})}$ of homogeneous linear recurrences (\ref{Gamk}) we investigate is (\ref{charpoly}). Assume that it has decomposition 
\begin{equation}\label{charpoly1}
	k_\M(x)=(x-\alpha_1)^{m_1}(x-\alpha_2)^{m_2}\cdots (x-\alpha_h)^{m_h}
\end{equation}
into linear factors, where $\alpha_i\in\C$ ($i=1,2,\dots,h$) are distinct, and $m_1+m_2+\cdots+m_h=k$. For each pair $(i,j)$, where $i\in\{1,2,\dots,h\}$ and 
$j\in\{0,1,\dots,m_i-1\}$ consider the sequences
\begin{equation}\label{seqs}
	S_n^{(i,j)}=n^j\alpha_i^n, \quad(n\ge0).
\end{equation}
Note that the range for $j$ depends on the value of $i$, but we do not indicate this fact in the notation. 

Sequences (\ref{seqs}) are belonging to $\Gamma_k^{(\gamma_0,\gamma_1,\dots,\gamma_{k-1})}$, this fact is based on the following two observations, or simply we may refer to the theory of difference equations, see e.g. the classical book of Milne-Thomson \cite[Chapter XIII]{M-T}. Firstly, if the multiplicity of an $\alpha_i$ is $m_i$, then $\alpha_i$ is also a root of the derivatives $k^{(j)}_\M(x)$ for $j=0,1,\dots,m_i-1$. Secondly, from the theory of Stirling number of second kind we know that $n^j=\sum_{v=0}^j{j\brace v}(n)_v$, where 
${j\brace v}$ denotes a Stirling number of second kind, and $(n)_v=n(n-1)\cdots(n-v+1)$ is a falling factorial.

Another important argument is that the sequences in (\ref{seqs}) are $\C$-linearly independent, consequently they form a basis of $\Gamma_k^{(\gamma_0,\gamma_1,\dots,\gamma_{k-1})}$. This follows again from   \cite[Chapter XIII]{M-T}. (A similar result is given in \cite{K}.)

Analogously to $\G$ we define the matrix $\G_S$ by the initial values of the sequences $S_n^{(i,j)}$ such that the first row contains the initial values of $S_n^{(1,0)}$, etc. The matrix $\G_S$ can be split into horizontal stripes such that stripe  $i\in\{1,2,\dots,h\}$ contains $m_i$ rows. The enumeration of stripes begins on the top of the matrix, and the entries of stripe $i$ are given as follows:
\begin{equation*}
	\G_S=\left[
	\begin{array}{cccccccc}
		&&&&&\vdots&& \\ \hdashline
		&&&&&\vdots&& \\ \hdashline
		1 & \alpha_i & \alpha_i^2 & \alpha_i^3 & \dots & & \dots & \alpha_i^{k-1} \\
		0 & \alpha_i & 2\alpha_i^2 & 3\alpha_i^3 & \dots & & \dots & (k-1)\alpha_i^{k-1} \\
		0 & \alpha_i & 2^2\alpha_i^2 & 3^2\alpha_i^3 & \dots & & \dots & (k-1)^2\alpha_i^{k-1} \\
		\vdots & \vdots & \vdots & \vdots & & \ddots & & \vdots \\
		0 & \alpha_i & 2^{m_i-1}\alpha_i^2 & 3^{m_i-1}\alpha_i^3 & \dots & & \dots & (k-1)^{m_i-1}\alpha_i^{k-1} \\ \hdashline
		&&&&&\vdots&& \\ \hdashline
		&&&&&\vdots&& \\ 
	\end{array}	
	\right].
\end{equation*}
Let $M_i=m_1+m_2+\cdots+m_{i-1}$. The numbering of rows of stripe $i$ are $M_i+(j+1)$, where $j=0,\dots,m_i-1$. The matrix $\G_S^\star$ is analogous to $\G^\star$, too.
We define even the block diagonal matrix 
$$
\B = \left[
\begin{array}{cccc}
	\B_1 & \O  & \cdots & \O \\
	\O & \B_2  & \cdots & \O \\
	\vdots & \vdots & \ddots & \vdots \\
	\O & \O & \cdots &  \B_h
\end{array}
\right]\in\C^{k\times k}\;{\rm with}
\;\B_i = \left[
\begin{array}{cccccc}
	\alpha_i & 0  & \cdots & 0 & \cdots & 0\\
	\alpha_i & \alpha_i  & \cdots & 0 & \cdots & 0\\
	\vdots & \vdots & \ddots & \vdots & \vdots & \vdots\\
	\alpha_i & \cdots  & \cdots & \alpha_i & \cdots & 0\\
	\vdots & \vdots & b_{u,v}^{(i)} & \vdots & \ddots & \vdots\\	
	\alpha_i & \cdots & \cdots & \cdots & \cdots & \alpha_i
\end{array}
\right]\in\C^{m_i\times m_i}
$$
for $i=1,2,\dots, h$
such that the general term $b_{u,v}^{(i)}=\binom{u}{v}\alpha_i$ is the entry of the $(u+1)$th row and $(v+1)$th column of $\B_i$ with the conditions $0\le v\le u \le m_i-1$. Here the minors $\O$ are suitable zero matrices. Since $\B_i$ is a lower diagonal matrix then $\B$ is so. 

Using the notation above we will show that 
\begin{equation*}
	\G_S^*=\B\G_S.
\end{equation*}
Note that block $\B_i$ has influence only on stripe $i$ of $\G_S$ during carrying out the matrix multiplication $\B\G_S$.
In stripe $i$ we have exactly the rows $M_i+1+u$, where $u=0,1,\dots,m_i-1$. Our deal is to calculate the dot product\footnote{The dot product of the vectors $\bar{\mu}=[\mu_1,\ldots,\mu_k]$ and $\bar{\nu}^\top=[\nu_1,\ldots,\nu_k]^\top$ is the complex number $\sum_{j=1}^k \mu_j\nu_j$.}  of row $M_i+1+u$ from $\B$ and column $\ell$ from $\G_S$, where $\ell=1,2,\dots, k$. 
This is given in details by
\begin{equation*}
	\left[
	\begin{array}{cccccccccc}
		0 & \cdots & 0  & \binom{u}{0}\alpha_i & \binom{u}{1}\alpha_i & \cdots & \binom{u}{u}\alpha_i  & 0 & \cdots & 0\\
	\end{array}
	\right]\cdot
	\left[
	\begin{array}{c}
		\vdots \\ \hdashline
		\alpha_i^{\ell-1} \\
		(\ell-1)\alpha_i^{\ell-1} \\
		\vdots \\
		(\ell-1)^u\alpha_i^{\ell-1} \\
		\vdots \\
		(\ell-1)^{m_i-1}\alpha_i^{\ell-1} \\ \hdashline
		\vdots
	\end{array}
	\right].
\end{equation*}
Hence the dot product equals
\begin{equation*}
	\binom{u}{0}\alpha_i^\ell+\binom{u}{1}(\ell-1)\alpha_i^\ell+\cdots+\binom{u}{u}(\ell-1)^u\alpha_i^\ell 
	=\left((l-1)+1\right)^u\alpha_i^\ell=\ell^u\alpha_i^\ell=S_\ell^{(i,u)}.
\end{equation*}
Recall that $S_\ell^{(i,u)}$ is the general term of the matrix $\G_S^*$ such that it is the element of the $(u+1)$th row of stripe $i$ in the column $\ell$ ($\ell=1,2,\dots, k$).

Subsequently, we have
$$
\M_S=\G_S^*\G_S^{-1}=\left(\B\G_S\right)\G_S^{-1}=\B,
$$
and then $\M_S^n=\B^n$ follows. We know that $\B$ is a lower diagonal matrix. Thus $\det(\B)$ is the product of the elements lying in the main diagonal, i.e.
\begin{equation*}
	\det(\B)=\alpha_1^{m_1}\alpha_2^{m_2}\cdots \alpha_h^{m_h}.
\end{equation*}
Finally,
\begin{equation*}
	\det\left(\M_S^n\right)=\left(\det(\B)\right)^n=\alpha_1^{m_1n}\alpha_2^{m_2n}\cdots \alpha_h^{m_hn}
	= \left(S_n^{(1,0)}\right)^{m_1}\left(S_n^{(2,0)}\right)^{m_2}\cdots \left(S_n^{(h,0)}\right)^{m_h}.
\end{equation*}
So we have found that $\det\left(\M_S^n\right)$
is a corresponding value of the decomposable form
$$
F_S(y_1,y_2,\dots,y_h)=y_1^{m_1}y_2^{m_2}\cdots y_h^{m_h}
$$ 
at the point
$\left((S_n^{(1,0)}),(S_n^{(2,0)}),\cdots, (S_n^{(h,0)})\right)$.
Thus the assertion of the theorem holds for this particular choice of the recurrences.
\bigskip

\noindent{\bf Part 2.}
Recall that sequences (\ref{kseq}) form a basis in $\Gamma_k^{(\gamma_0,\gamma_1,\dots,\gamma_{k-1})}$ because $\Delta\ne0$. Hence there uniquely exists a matrix $\A=[a_{i,j}]\in\C^{h\times k}$ such that

\begin{equation}\label{last}
	\left[
	\begin{array}{cccc}
		S_n^{(1,0)} & S_n^{(2,0)} & \cdots  & S_n^{(h,0)}
	\end{array}
	\right]^\top=\A \left[
	\begin{array}{cccc}
		G_n^{(1)} & G_n^{(2)} & \cdots  & G_n^{(k)}
	\end{array}
	\right]^\top.
\end{equation}
Consider now the substitution $\bar{y}=[y_1,y_2,\ldots,y_h]^\top = {\A} [x_1,x_2,\ldots,x_k]^T={\A}\bar{x}$. Denoting by $\bar{a}_i$ the $i$th row vector of $\A$ we compute the dot product
$$
y_i = \bar{a}_i[x_1,x_2,\ldots,x_k]^\top  
$$
for every $i=1,2,\ldots,h$. Now we can define $F(x_1,\ldots,x_k)=F(\bar{x})$ as follows:
$$
F_S(\bar{y})=y_1^{m_1}y_2^{m_2}\cdots y_h^{m_h}=\prod_{i=1}^h \left(\bar{a}_i[x_1,x_2\ldots,x_k]^\top\right)^{m_i}=F(\bar{x}).
$$ 
Obviously, $F$ is a decomposable form of degree $k$ with the property
$$
F(G_n^{(1)},\ldots,G_n^{(k)}) = 
((-1)^{k} k_{\M}(0))^n=\delta^n,\quad(n=0,1,\ldots).
$$
Thus the proof is complete.	$\Box$

\medskip

\begin{remark}
If every zero of the characteristic polynomial $k_\M(x)$ is simple (i.e.~$h=k$), then the matrix ${\G}_S$ is simplified to
\begin{equation*}
	{\G}_S = \left[
	\begin{array}{cccc}
		1 & \alpha_1  & \cdots &  \alpha_1^{k-1} \\
		1 & \alpha_2  & \cdots &  \alpha_2^{k-1} \\
		\vdots & \vdots & \ddots & \vdots \\
		1 & \alpha_k  & \cdots &  \alpha_k^{k-1} 
	\end{array}
	\right]=\A\G
\end{equation*}
(the second equality is implied by (\ref{last})). Consequently, $\A=\G_S\G^{-1}$, and then 
\begin{equation*}
  \bar{y}=\G_S\G^{-1}\bar{x},
\end{equation*}	
Finally, we obtain $F_S(\bar{y})=F_S(\G_S\G^{-1}\bar{x})=F(\bar{x})$.
\end{remark}

At the end of this section, we show a general example to demonstrate the power of our results.
\bigskip

\noindent{\bf Example.} Consider the sequences of (\ref{kseq}) such that the characteristic polynomial has factorization (\ref{charpoly1}). Assume that we have the specific initial vectors
\begin{equation}\label{defg_new}
	\bg_0=\left[\begin{array}{c}
		1 \\ 0 \\ \vdots \\ 0 \\
	\end{array}
	\right],\;
	\bg_1=\left[\begin{array}{c}
		0 \\ 1 \\ \vdots \\ 0 \\
	\end{array}
	\right],\dots,\;
	\bg_{k-1}=\left[\begin{array}{c}
		0 \\ 0  \\ \vdots \\ 1 \\
	\end{array}
	\right].
\end{equation}
Using the advantageous properties of this orthonormal basis, we can easily show that
$$
S_n^{(j,0)}=1\cdot G_n^{(1)} + \alpha_j\cdot G_n^{(2)} + \cdots + \alpha_j^{k-1}\cdot G_n^{(k)}, \quad(j=1,2,\dots,h).
$$
Consequently, 
$$
y_j=x_1+ \alpha_jx_2 + \cdots + \alpha_j^{k-1}x_k, \quad(j=1,2,\dots,h),
$$
and then
\begin{equation*}
	{F}(x_1,x_2,\dots,x_k)=\prod_{i=1}^h(x_1+ \alpha_jx_2 + \cdots + \alpha_j^{k-1}x_k)^{m_j}
\end{equation*}
follows. Hence the complete factorization of the form ${F}$ is determined. Note that the simplicity depends on the features of the orthonormal basis (\ref{defg_new}). 

\section{Proof of the corollaries and of Theorem \ref{t:indep}} \label{s:corrolaries}

{\it Proof of Corollary \ref{c:ugyanaz}.} The entries of the vectors $\bar{g}_j$, $j=0,\dots,k-1$ are the initial values of the sequences $(G_{n+j})$, $j=0,\dots, k-1$. At the beginning of the proof of Theorem \ref{t:fo} we pointed out that the  linear independence of the vectors $\bar{g}_j$ and the linear independence of the sequences $(G_{n+j})$ is equivalent. Thus  the sequences $(G_{n+j})$, $j=0,\dots, k-1$ are linearly independent and Theorem \ref{t:fo} implies the existence of a homogeneous form $F$, while Theorem \ref{t:decompose} justifies the decomposability of $F$.
$\Box$
\medskip

{\it Proof of Corollary \ref{c:indep}.} If the sequences  $(G_n^{(j)})$, $j=1,\dots, k$ are linearly dependent, then enlarging their set with the sequence $(\gamma_0^n)$ the linear dependent property remains true, and  this is a very special algebraic dependence. Otherwise, if $(G_n^{(j)}), j=1,\dots, k$ are linearly independent, then Theorem \ref{t:fo}, more precisely equation \eqref{mainresult} verifies the algebraic dependence. $\Box$
\medskip

{\it Proof of Theorem \ref{t:indep}.} Assume that $P(X)=z_kX^k + z_{k-1}X^{k-1}+\ldots +z_0$, where $z_0\not= 0$. The differential equation with characteristic polynomial $P(X)$ has the form
\begin{equation} \label{e:diffe}
z_kf(z)^{(k)} + z_{k-1}f(z)^{(k-1)}+\ldots +z_0 = 0.
\end{equation}
Here $f^{(\ell)}$ denotes the $\ell$th derivative of $f$. The pair-wisely distinct zeros of $P(X)$ are assigned by $\alpha_1,\alpha_2,\dots,\alpha_k$. The function $f(z)$ is a solution of \eqref{e:diffe} if and only if it has the form 
$$
f(z) = f_1\alpha_1^z + f_2\alpha_2^z +\ldots + f_k\alpha_k^z,
$$
with complex numbers $f_1,f_2,\dots,f_k$.

The set of solutions of \eqref{e:diffe} form a $\mathbb{C}$-vector space with respect to the addition of functions and multiplication with complex numbers. We denote it by $S(P)$. The functions $\alpha_1^z,\alpha_2^z,\dots,\alpha_k^z$ are linearly independent, thus the dimension of $S(P)$ is $k$. 

Assume that $G_j(z), j\le k_0$ are pairwise different solutions of \eqref{e:diffe}. If they are $\mathbb{C}$-linearly dependent, which holds always when $k_0>k$, then the same is true if we enlarge their set with the function $((-1)^kP(0))^z$. It remains to examine the case when $k_0=k$ and the functions $G_j(z), j\le k_0$ are linearly independent. Then they form a bases of $S(P)$. The functions $\alpha_1^z,\alpha_2^z,\dots,\alpha_k^z$ belong to $S(P)$, thus there exist complex numbers $a_{ij}, 1\le i,j\le k$ such that
\begin{equation}\label{qqriq}
\alpha_i^z = \sum_{j=1}^k a_{ij} G_j(z).
\end{equation}

On the other hand 
$$
\alpha_1^z\alpha_2^z\cdots\alpha_k^z = (\alpha_1\alpha_2\cdots\alpha_k)^z = ((-1)^k P(0))^z.
$$
Inserting here the linear relations (\ref{qqriq}) we obtain the statement. $\Box$

\section{Binary and ternary recurrences} \label{s:negy}

\subsection{Binary recurrences, a general identity}

This subsection is devoted to present the general identity (\ref{fineagain}) below with two binary recurrences satisfying the same recurrence relation. This is a special case of (\ref{mainresult}). On the other hand it provides a common generalization of (\ref{Gen1})-(\ref{jan13}).
At the end of this subsection a few examples will illustrate (\ref{fineagain}). In the computations, we will follow the arguments of the previous sections.

Assume $k=2$. For simplicity suppose that the two recurrent sequences are $(G_n)$ and $(H_n)$, their initial values are $G_0,G_1$ and $H_0,H_1$, respectively, furthermore put $\gamma_1=A$, $\gamma_0=B$. Then the two sequences above belong to
\begin{equation*}\label{Gam2}
	\Gamma_2^{(B,A)}=\left\{(x_n)_{n\in\Z}\in\C^\infty\mid x_n=Ax_{n-1}+Bx_{n-2},\;n\ge 2\right\}.
\end{equation*}
Now
\begin{equation*}
	\G=\left[\begin{array}{cc}
		G_0 & G_1 \\
		H_0 & H_1
	\end{array}
	\right],\quad\quad \G^\star=\left[\begin{array}{cc}
		G_1 & AG_1+BG_0 \\
		H_1 & AH_1+BH_0
	\end{array}
	\right],
\end{equation*}
and assume that $\Delta=\det(\G)=G_0H_1-G_1H_0\ne0$. Then
\begin{equation*}
	\G^{-1}=\frac{1}{\Delta}\left[\begin{array}{cc}
		{H_1} & -{G_1} \\
		-{H_0}  & {G_0}
	\end{array}
	\right].
\end{equation*}
Introduce the notation $C_G=G_1^2-AG_0G_1-BG_0^2$ (see also before the equations \eqref{Gen1}-\eqref{jan13}), and analogously $C_H=H_1^2-AH_0H_1-BH_0^2$. We even define $E_{10}=G_1H_1-AG_1H_0-BG_0H_0$, and $E_{01}=G_1H_1-AG_0H_1-BG_0H_0$. One can easily see that
\begin{equation}\label{Mexpl}
	\M=\G^\star\G^{-1}=
	\frac{1}{\Delta}\left[\begin{array}{cc}
		{E_{10}} & {-C_G}  \\
		{C_H} & {-E_{01}}
	\end{array}
	\right].
\end{equation}
Clearly, $k_\M(x)=x^2-Ax-B=(x-\alpha)(x-\beta)$, $\det(\M)=-B$. Obviously, $\det(\M^n)=(-B)^n$, $A=\alpha+\beta$ and $B=-\alpha\beta$. (Note that $\alpha$ and $\beta$ are not necessarily distinct.)

The element sequences of the powers of matrix $\M$ satisfy
\begin{equation*}
	m_{i,j}^{(n)}=Am_{i,j}^{(n-1)}+Bm_{i,j}^{(n-2)},\qquad (n\ge 2;\;1\le i,j\le2),
\end{equation*}
the initial values are clear from $\M^0=\mathbf{I}$ and from $\M$ (here $\mathbf{I}$ is the $2\times2$ unit matrix). Now
\begin{equation} \label{sys}
	\left[\begin{array}{c}
		{m_{i,j}^{(0)}} \\
		{m_{i,j}^{(1)}}
	\end{array}
	\right]=
	c_{i,j}^{(1)}
	\left[\begin{array}{c}
		G_0 \\
		G_1
	\end{array}
	\right]+
	c_{i,j}^{(2)}
	\left[\begin{array}{c}
		H_0 \\
		H_1
	\end{array}
	\right],\quad\quad(1\le i,j\le2).
\end{equation}

Once we have the solutions to (\ref{sys}) in $c_{i,j}^{(1)}$ and $c_{i,j}^{(2)}$, then $m_{i,j}^{(n)}=c_{i,j}^{(1)}G_n+c_{i,j}^{(2)}H_n$ holds. A straightforward calculation shows that
\begin{equation*}
	\mathbf{M}^n=\frac{1}{\Delta}\left[\begin{array}{lr}
		(-H_0m_{1,1}^{(1)}+H_1)G_n+(G_0m_{1,1}^{(1)}-G_1)H_n & -H_0m_{1,2}^{(1)}G_n+G_0m_{1,2}^{(1)}H_n  \\
		-H_0m_{2,1}^{(1)}G_n+G_0m_{2,1}^{(1)}H_n  & (-H_0m_{2,2}^{(1)}+H_1)G_n+(G_0m_{2,2}^{(1)}-G_1)H_n
	\end{array}
	\right],
\end{equation*}
where we used the values $m_{1,1}^{(0)}=m_{2,2}^{(0)}=1$, $m_{1,2}^{(0)}=m_{2,1}^{(0)}=0$. In fact, to determine $\det(\mathbf{M}^n)$ we do not need the exact values $m_{i,j}^{(1)}$ given in (\ref{Mexpl}), only the identities $m_{1,1}^{(1)}m_{2,2}^{(1)}-m_{1,2}^{(1)}m_{2,1}^{(1)}=\det(\M)=-B$ and
$m_{1,1}^{(1)}+m_{2,2}^{(1)}={\rm tr}(\M)=A$.
Indeed, if we figure out simply the determinant of the matrix $\mathbf{M}^n$, and collect the coefficient of the terms $G_n^2$, $G_nH_n$, and $H_n^2$, respectively, then
the key moment of the simplification is the application of these identities.  Finally, the computations lead to
\begin{equation*}
	(-B)^n=\det\left(\mathbf{M}^n\right)=\frac{C_H}{\Delta^2}G_n^2+\frac{C_{GH}}{\Delta^2}G_nH_n+\frac{C_G}{\Delta^2}H_n^2,
\end{equation*}
where $C_{GH}=-(E_{10}+E_{01})$. Moreover one can show that $C_{GH}$ can be given by using the corresponding associate sequences: $C_{GH}=G_0\widehat{H}_1-G_1\widehat{H}_0=H_0\widehat{G}_1-H_1\widehat{G}_0$.
Then we have the desired equality given in
\begin{theorem}\label{main2}
	The terms of the recurrences above satisfy
	\begin{equation} \label{fineagain}
		C_HG_n^2+C_{GH}G_nH_n+C_GH_n^2=(-B)^n\Delta^2.
	\end{equation}
\end{theorem}

This is a nice common generalization of (\ref{Gen1})-(\ref{jan13}). 
Indeed, observe that Theorem \ref{main2} leads to (\ref{Gen1}) whenever $(H)$ is the associate sequence of $(G)$. Really, it is easy to show that $C_{\widehat{G}}=-DC_G$, $C_{G\widehat{G}}=0$ (i.e.~$C_{GH}$ vanishes if $(H)=(\widehat{G})$), and in this particular case $\Delta=-2C_G$ holds. Insert these values into (\ref{fineagain}) to get immediately (\ref{Gen1}). For (\ref{jan13}) let $H_n=G_{n+1}$, the details here are omitted. 

The binary quadratic form on left-hand side of (\ref{fineagain}) is trivially decomposable. The decomposition is formulated by
\begin{equation*}
	\left((H_1-\alpha H_0)G_n-(G_1-\alpha G_0)H_n\right))\cdot\left((H_1-\beta H_0)G_n-(G_1-\beta G_0)H_n\right)=(-B)^n\Delta^2.
\end{equation*}
We now illustrate Theorem {\ref{main2} by five examples given in Table \ref{tab1}.	

\begin{table}[h]\label{tab1}
\renewcommand{\arraystretch}{1.2}
\begin{tabular}{|c|c|c|c|} \hline
	
	$(A,B)$  & $\mathbf{G}$ & $\mathbf{M}$ & $(-B)^n\Delta^2=C_HG_n^2+C_{GH}G_nH_n+C_GH_n^2$  \\ \hline\hline
	
	$(0,4)$  &  $\left[\begin{array}{cc}
		1 & 2 \\
		2 & 3
	\end{array}\right]$   & $\left[\begin{array}{cc}
		2 & 0 \\
		7 & -2
	\end{array}\right]$  & $(-4)^n(-1)^2=-7G_n^2+4G_nH_n=G_n(4H_n-7G_n)$  \\ \hline

	$(2,-1)$ & $\left[\begin{array}{cc}
		2 & 3 \\
		4 & 5
	\end{array}\right]$  & $\left[\begin{array}{cc}
		1/2 & 1/2 \\
		-1/2 & 3/2
	\end{array}\right]$  & $1^n(-2)^2=G_n^2-2G_nH_n+H_n^2=(G_n-H_n)^2$  \\ \hline
	
	$(7,-10)$   & $\left[\begin{array}{cc}
		0 & 1 \\
		2 & 7
	\end{array}\right]$  & $\left[\begin{array}{cc}
		7/2 & 1/2 \\
		9/2 & 7/2
	\end{array}\right]$  & $10^n(-2)^2=-9G_n^2+H_n^2$  \\ \hline

$(7,-10)$  & $\left[\begin{array}{cc}
	1 & 2 \\
	1 & 5
\end{array}\right]$ & $\left[\begin{array}{cc}
	2 & 0 \\
	0 & 5
\end{array}\right]$  & $10^n3^2=9G_nH_n$    \\ \hline

	$(4,-1)$  & $\left[\begin{array}{cc}
	1 & 2 \\
	3 & 4
\end{array}\right]$   & $\left[\begin{array}{cc}
	13/2 & -3/2 \\
	23/2 & -5/2
\end{array}\right]$  & $1^n(-2)^2=-23G_n^2+18G_nH_n-3H_n^2$  \\ \hline
	
\end{tabular}
\caption{\label{tab1}Binary recurrence examples.}
\end{table}
One or two coefficients from $\{C_G, C_{H}, C_{GH}\}$ may vanish in (\ref{fineagain}), which provides a large variety of identities.

\bigskip

\subsection{Ternary recurrences}

 Instead of pointing on the complicated general formula it is better to present an example which is typical of the ternary forms we are studying. Examples for higher order linear recurrences can be analogously produced.
 
 There are only a few papers in the literature which work with three different ternary recurrences satisfying the same recurrence rule. An example is \cite{S}, which studies Narayana sequence $(A_n)$, Narayana-Lucas sequence $(B_n)$, and Narayana-Perrin sequence $(C_n)$.
These recurrences satisfy 
\begin{equation*}
	x_n=x_{n-1}+x_{n-3}
\end{equation*}
with initial values given by
\begin{equation*}
	A_0=0,A_1=1,A_2=1;\quad
	B_0=3,B_1=1,B_2=1;\quad
	C_0=3,C_1=0,C_2=2,
\end{equation*}
respectively. Let $\alpha_i$ ($i=1,2,3$) denote the simple zeros of the characteristic polynomial $k(x)=x^3-x^2-1$ such that $\alpha_1\in\R$ and $\alpha_3$ is the complex conjugate of $\alpha_2$.
Following the the method we described in details it leads to the diophantine equation
\begin{equation}\label{tiz}
	-187x^3+159x^2y-45x^2z-189xy^2+306xyz-117xz^2+y^3-45y^2z+63yz^2-27z^3 = -216
\end{equation}
possessing infinitely many integer solutions $x=A_n$, $y=B_n$, $z=C_n$.
In this case, $\Delta=-6$, and the base of the exponential term is $1$. Note that all the 10 coefficients in the corresponding $\tilde{F}$ are non-zero.
The ternary form of (\ref{tiz}) can be decomposed in the algebraic number field $\Q[\alpha_1,\alpha_2]$ as 
\begin{equation*}
	-\prod_{i=1}^{3}\left((3\alpha_i^2+3\alpha_i-2)x+(-3\alpha_i^2+3\alpha_i+2)y+(3\alpha_i^2-3\alpha_i)z\right).
\end{equation*}

\subsection*{Acknowledgments}	
We thank Robert Tichy for his help in the formulation of Theorem 3.
	K. Gy\H{o}ry and L. Szalay were supported by the Hungarian National Foundation for Scientific Research Grant No.~128088, and No.~130909.
	For L.~Szalay this research was supported by the Slovak Scientific Grant Agency VEGA 1/0776/21. 

\end{document}